\documentclass{amsart}
\usepackage{latexsym,amsmath,amssymb,amsthm}

\theoremstyle{definition}
\theoremstyle{remark}

\numberwithin{equation}{section}

\begin{document}

\title[WEIGHTED CONDITIONAL TYPE OPERATORS ON ORLICZ SPACES]
{WEIGHTED CONDITIONAL TYPE OPERATORS BETWEEN DIFFERENT ORLICZ SPACES}

\author{\bf Y. Estaremi}
\address{ Y. Estaremi}
\address{ Department of Mathematics, Payame Noor University (PNU), P. O. Box: 19395-3697, Tehran- Iran}
 \email{estaremi@gmail.com}



\begin{abstract}
In this note we consider weighted conditional type operators between different Orlicz spaces and generalized conditional type H\"{o}lder inequality that we defined in \cite{e3}. Then we give some necessary and sufficient conditions for boundedness of weighted conditional type operators. As a consequence we characterize boundedness of weighted conditional type operators and multiplication operators between different $L^p$-spaces. Finally, we give some upper and lower bounds for essential norm of weighted conditional type operators.
\end{abstract}

\maketitle

\section{\sc\bf Introduction and Preliminaries}
  The continuous convex function
$\Phi:\mathbb{R}\rightarrow\mathbb{R}$ is called a Young's function whenever

(1)$\Phi(x)=0$ if and only if $x=0$.

(2) $\Phi(x)=\Phi(-x)$.

(3) $\lim_{x\rightarrow\infty}\frac{\Phi(x)}{x}=\infty$, $\lim_{x\rightarrow\infty}\Phi(x)=\infty$.

 With each Young's function $\Phi$ one can associate another
 convex function $\Phi^{*}:\mathbb{R}\rightarrow\mathbb{R}^{+}$ having similar properties, which is defined by $$\Phi^{*}(y)=\sup\{x|y|-\Phi(x):x\geq0\}, \ \ y\in\mathbb{R}.$$
Then $\Phi^{\ast}$ is called the complementary Young's function of $\Phi$. The following properties also are immediate from the definition of Young functions.\\

\vspace*{0.3cm} {\bf Proposition 1.1.}\cite{raor} Let $\Phi$ be a Young's function. Then $\Phi$, $\Phi^{*}$ are strictly increasingly so that their inverses $\Phi^{-1}$, $\Phi^{\ast^{-1}}$ are uniquely defined and\\
(i) $\Phi(a)+\Phi(b)\leq \Phi(a+b)$, $\Phi^{-1}(a)+\Phi^{-1}(b)\geq \Phi^{-1}(a+b)$, $a,b\geq0$,\\
(ii) $a<\Phi^{-1}(a)\Phi^{*^{-1}}(a)\leq 2a$, $a\geq0$.\\

A Young's function $\Phi$ is said to satisfy the
$\bigtriangleup_{2}$ condition (globally) if $\Phi(2x)\leq
k\Phi(x), \ x\geq x_{0}\geq0 \ \  (x_{0}=0)$ for some constant
$k>0$. Also, $\Phi$ is said to satisfy the
$\bigtriangleup'(\bigtriangledown')$ condition, if $\exists c>0$
$(b>0)$ such that
$$\Phi(xy)\leq c\Phi(x)\Phi(y), \ \ \ x,y\geq x_{0}\geq 0$$
$$(\Phi(bxy)\geq \Phi(x)\Phi(y), \ \ \ x,y\geq y_{0}\geq 0).$$
If $x_{0}=0(y_{0}=0)$, then these conditions are said to hold
globally. If $\Phi\in \bigtriangleup'$, then $\Phi\in
\bigtriangleup_{2}$.

 Let $\Phi, \Psi$ be two Young's
functions, then $\Phi$ is stronger than $\Psi$,
$\Phi\succ\Psi$ [or $\Psi\prec\Phi$] in symbols, if
$$\Psi(x)\leq\Phi(ax), \ \ \ x\geq x_0\geq0$$
for some $a\geq0$ and $x_0$, if $x_0=0$ then this condition is
said to hold globally. A detailed discussion and verification of these properties may be found in \cite{raor}.

Let $(\Omega, \Sigma, \mu)$ be a measure space and $\Phi$ be a Young's function, then the set of $\Sigma$-measurable functions
 $$L^{\Phi}(\Sigma)=\{f:\Omega\rightarrow \mathbb{C}:\exists k>0, \int_{\Omega}\Phi(k|f|)d\mu<\infty\}$$
 is a Banach space, with respect to the norm $N_{\Phi}(f)=\inf\{k>0:\int_{\Omega}\Phi(\frac{f}{k})d\mu\leq1\}$.
 $(L^{\Phi}(\Sigma), N_{\Phi}(.))$ is called Orlicz space. If $\Phi\in \bigtriangleup_{2}$, then the dual space of $L^{\Phi}(\Sigma)$ is equal to $L^{\Phi^*}(\Sigma)$. The
 usual convergence in the orlicz space $L^{\Phi}(\Sigma)$ can be
 introduced in term of the orlicz norm $N_{\Phi}(.)$ as
 $u_n\rightarrow u$ in $L^{\Phi}(\Sigma)$ means $N_{\Phi}(u_n-u)\rightarrow
 0$. Also, a sequence $\{u_n\}^{\infty}_{n=1}$ in
 $L^{\Phi}(\Sigma)$ is said to converges in $\Phi$-mean to $u\in
 L^{\Phi}(\Sigma)$, if
  $$\lim_{n\rightarrow\infty}
 I_{\Phi}(u_n-u)=\lim_{n\rightarrow\infty}\int_{\Omega}\Phi(|u_n-u|)d\mu=0.$$

For a sub-$\sigma$-finite algebra $\mathcal{A}\subseteq\Sigma$, the
conditional expectation operator associated with $\mathcal{A}$ is
the mapping $f\rightarrow E^{\mathcal{A}}f$, defined for all
non-negative, measurable function $f$ as well as for all $f\in
L^1(\Sigma)$ and $f\in L^{\infty}(\Sigma)$, where
$E^{\mathcal{A}}f$, by the Radon-Nikodym theorem, is the unique
$\mathcal{A}$-measurable function satisfying
$$\int_{A}fd\mu=\int_{A}E^{\mathcal{A}}fd\mu, \ \ \ \forall A\in \mathcal{A} .$$
As an operator on $L^{1}({\Sigma})$ and $L^{\infty}(\Sigma)$,
$E^{\mathcal{A}}$ is idempotent and
$E^{\mathcal{A}}(L^{\infty}(\Sigma))=L^{\infty}(\mathcal{A})$ and
$E^{\mathcal{A}}(L^1(\Sigma))=L^1(\mathcal{A})$. Thus it can be
defined on all interpolation spaces of $L^1$ and $L^{\infty}$ such
as, Orlicz spaces \cite{besh}. We say the measurable function $f$ is conditionable with respect to $\sigma$-subalgebra $\mathcal{A}\subseteq \Sigma$ if
$E^{\mathcal{A}}(f)$ is defined. If there is no possibility of
confusion, we write $E(f)$ in place of $E^{\mathcal{A}}(f)$. This
operator will play a major role in our work and we list here some
of its useful properties:

\vspace*{0.2cm} \noindent $\bullet$ \  If $g$ is
$\mathcal{A}$-measurable, then $E(fg)=E(f)g$.

\noindent $\bullet$ \ $\varphi(E(f))\leq E(\varphi(f))$, where
$\varphi$ is a convex function.

\noindent $\bullet$ \ If $f\geq 0$, then $E(f)\geq 0$; if $f>0$,
then $E(f)>0$.

\noindent $\bullet$ \ For each $f\geq 0$, $S(f)\subseteq S(E(f))$,
where  $S(f)=\{x\in X; f(x)\neq 0\}$.

\vspace*{0.2cm}\noindent A detailed discussion and verification of
most of these properties may be found in \cite{rao}. We recall
that an $\mathcal{A}$-atom of the measure $\mu$ is an element
$A\in\mathcal{A}$ with $\mu(A)>0$ such that for each
$F\in\mathcal{A}$, if $F\subseteq A$, then either $\mu(F)=0$ or
$\mu(F)=\mu(A)$. A measure space $(\Omega,\Sigma,\mu)$ with no atoms is
called a non-atomic measure space. It is well-known fact that every
$\sigma$-finite measure space $(\Omega, \Sigma,\mu)$ can be
partitioned uniquely as $\Omega=\left
(\bigcup_{n\in\mathbb{N}}C_n\right )\cup B$, where
$\{C_n\}_{n\in\mathbb{N}}$ is a countable collection of pairwise
disjoint $\Sigma$-atoms and $B$, being disjoint from each $C_n$,
is non-atomic \cite{z}.\ \ \

Let $f\in L^{\Phi}(\Sigma)$.  It is not difficult to see that $\Phi(E(f))\leq
E(\Phi(f))$ and so by some elementary computations we get that $N_{\Phi}(E(f))\leq N_{\Phi}(f)$ i.e, $E$ is a
contraction on the Orlicz spaces.
As we defined in \cite{e3}, we say that the pair $(E, \Phi)$ satisfies the generalized conditional-type
H\"{o}lder-inequality (or briefly GCH-inequality) if there exists some positive constant $C$
such that for all $f\in L^{\Phi}(\Omega, \Sigma, \mu)$ and $g\in
L^{\Psi}(\Omega, \Sigma, \mu)$ we have
$$E(|fg|)\leq C \Phi^{-1}(E(\Phi(|f|)))\Phi^{*^{-1}}(E(\Phi^{*}(|g|))),$$
where $\Psi$ is the complementary Young's function of $\Phi$. There are many examples of the pair $(E, \Phi)$ that satisfy GCH-inequality in \cite{e3}.\\
This work is the continuance of \cite{e3}. In this paper we investigate boundedness of weighted conditional type operators between different Orlicz spaces by considering GCH-inequality. The results of the section 2 generalizes some results of \cite{e3} and \cite{taky}. In section 3 we find some upper and lower bounds for weighted conditional type operators on Orlicz spaces.

\section{ \sc\bf Bounded weighted conditional type operators }
 First we give a definition of weighted conditional type operator.

\vspace*{0.4cm} {\bf Definition 2.1.}
Let $(\Omega,\Sigma,\mu)$ be a $\sigma$-finite measure space and let $\mathcal{A}$ be a
$\sigma$-subalgebra of $\Sigma$ such that $(\Omega,\mathcal{A},\mathcal{A})$ is also $\sigma$-finite. Let $E$ be the corresponding conditional
expectation operator relative to $\mathcal{A}$. If $u \in L^0(\Sigma)$ (the spaces of $\Sigma$-measurable functions on $\Omega$) such that $uf$ is conditionable and $E(uf)\in L^{\Psi}(\Sigma)$ for all $f\in \mathcal{D}\subseteq L^{\Phi}(\Sigma)$, where $\mathcal{D}$ is a linear subspace, then the corresponding weighted conditional type operator (or WCT operator) is the linear transformation $R_{u}:\mathcal{D}\rightarrow L^{\Psi}(\Sigma)$ defined by $f\rightarrow E(uf)$.\\
In the first theorem we give some necessary conditions for boundedness of  $R_{u}:L^{\Phi}(\Sigma)\rightarrow L^{\Psi}(\Sigma)$, when $\Phi\preceq \Psi$  and some sufficient conditions, when $\Phi\succeq \Psi$.\\

\vspace*{0.3cm} {\bf Theorem 2.2.} Let WCT operator $R_{u}:\mathcal{D}\subseteq L^{\Phi}(\Sigma)\rightarrow L^{\Psi}(\Sigma)$ be well defined, then the followings hold.\\

(a) Let $\mu(\Omega)<\infty$ ($\mu(\Omega)=\infty$) and $\Phi\preceq \Psi$ (globally). Then

 (i) If $R_u$ is bounded from $L^{\Phi}(\Sigma)$ into $L^{\Psi}(\Sigma)$, then $E(u)\in L^{\infty}(\mathcal{A})$.

 (ii) If $\Psi\in \bigtriangleup'$(globally) and $R_u$ is bounded from $L^{\Phi}(\Sigma)$ into $L^{\Psi}(\Sigma)$, then  $\Psi^{*^{-1}}(E(\Psi^*(u)))\in L^{\infty}(\mathcal{A})$.

 (b) Let $\mu(\Omega)<\infty$ ($\mu(\Omega)=\infty$) and $\Psi\preceq \Phi$ (globally). Moreover, if $(E, \Psi)$ satisfies the GCH-inequality and $\Psi^{*^{-1}}(E(\Psi^*(u)))\in L^{\infty}(\mathcal{A})$, then $R_u$ is bounded. In this case, $\|R_u\|\leq C\|\Psi^{*^{-1}}(E(\Psi^*(u)))\|_{\infty}$, where the constant C comes from GCH-inequality.

\vspace*{0.3cm} {\bf Proof.} (a)-(i) Suppose that $E(u)\notin L^{\infty}(\mathcal{A})$. If we set $E_{n}=\{w\in \Omega:|E(u)(w)|>n\}$,  for all $n\in \mathbb{N}$, then $E_{n}\in \mathcal{A}$ and $\mu(E_{n})>0$.  Since $(\Omega, \mathcal{A}, \mu)$  has the finite subset property, we can assume that $0<\mu(E_{n})<\infty$, for all $n\in \mathbb{N}$.  By definition of $E_{n}$ we have  $$R_u(\chi_{E_{n}})=E(u\chi_{E_{n}})=E(u)\chi_{E_{n}}>n\chi_{E_{n}}.$$
 Since $\Phi\preceq \Psi$ and the Orlicz's norm is monotone, thus there exists a positive constant $c$ such that  $$\|R_u(\chi_{E_{n}})\|_{\Psi}\geq \frac{1}{c}\|R_u(\chi_{E_{n}})\|_{\Phi}>\frac{1}{c}\|n\chi_{E_{n}}\|_{\Phi}=\frac{n}{c}\|\chi_{E_{n}}\|_{\Phi}.$$
 This implies that $R_u$ isn't bounded. Therefor $E(u)$ should be essentially bounded.

 (a)-(ii) If $\Psi^{*^{-1}}(E(\Psi^*(u)))\notin L^{\infty}(\mathcal{A})$, then $\mu(E_{n})>0$, where $$E_{n}=\{w\in \Omega:\Psi^{*^{-1}}(E(\Psi^*(u)))(w)>n\}$$
and so $E_{n}\in \mathcal{A}$. Since $\Psi\in \bigtriangleup'$, then $\Psi^{*}\in\bigtriangledown'$, i.e., $\exists b>0$ such that $$\Psi^*(bxy)\geq\Psi^*(x)\Psi^*(y), \ \ \ x, y\geq0.$$
  Also, $\Phi, \Psi\in \triangle_{2}$. Thus $(L^{\Phi})^{\ast}=L^{\Phi^*}$ and $(L^{\Psi})^{\ast}=L^{\Psi^*}$ and so $T^{\ast}=M_{\bar{u}}:L^{\Psi^*}(\mathcal{A})\rightarrow L^{\Phi^*}(\Sigma)$, is also bounded. Hence for each $k>0$ we have
  \begin{align*}
  \int_{\Omega}\Psi^*(\frac{ku\chi_{E_{n}}}{N_{\Psi^*}(\chi_{E_{n}})})d\mu&=\int_{\Omega}\Psi^*(ku\chi_{E_{n}}\Psi^{*^{-1}}(\frac{1}{\mu(E_{n})}))d\mu\\
  &\geq\int_{E_{n}}\Psi^*(u)\Psi^*(\frac{ck\Psi^{*^{-1}}(\frac{1}{\mu(E_{n})})}{b})d\mu\\
  &\geq\left(\int_{E_{n}}E(\Psi^*(u))d\mu\right)\Psi^*(\frac{ck}{b^2})\Psi^*(\Psi^{*^{-1}}(\frac{1}{\mu(E_{n})}))\\
  &\geq\Psi^*(n)\mu(E_{n})\frac{1}{\mu(E_{n})}\Psi^*(\frac{ck}{b^2})\\
  &=\Psi^*(n)\Psi^*(\frac{ck}{b^2}).
  \end{align*}
Thus
$$\int_{\Omega}\Psi^*(\frac{ku\chi_{E_{n}}}{N_{\Psi^*}(\chi_{E_{n}})})d\mu=\int_{\Omega}\Psi^*(kM_{u}(f_{n}))d\mu
\geq\Psi^*(n)\Psi^*(\frac{k}{b^2})\rightarrow\infty$$ as $n\rightarrow\infty$, where $f_{n}=\frac{\chi_{E_{n}}}{N_{\Psi^*}(\chi_{E_{n}})}$. Thus
$N_{\Psi^*}(M_{u}(f_{n}))\rightarrow\infty$, as $n\rightarrow\infty$. Since $N_{\Phi^*}(M_{u}(f_{n}))\leq N_{\Psi^*}(M_{u}(f_{n}))$, then $N_{\Phi^*}(M_{u}(f_{n}))\rightarrow\infty$, as $n\rightarrow\infty$. This is a contradiction, since $M_{u}$ is bounded.

(b) Put $M=\|\Psi^{*^{-1}}(E(\Psi^*(u)))\|_{\infty}$. For $f\in L^{\Phi}(\Sigma)$ we have
\begin{align*}
\int_{\Omega}\Psi(\frac{E(uf)}{CMN_{\Phi}(f)})d\mu&=\int_{\Omega}\Psi(\frac{E(u\frac{f}{N_{\Phi}(f)})}{CM})d\mu\\
&\leq\int_{\Omega}\Psi(\frac{C \Psi^{*^{-1}}(E(\Psi(\frac{f}{N_{\Phi}(f)})))\Psi^{*^{-1}}(E(\Psi^*(u)))}{CM})d\mu\\
&\leq\int_{\Omega}\Psi(\Psi^{-1}(E(\Psi(\frac{f}{N_{\Phi}(f)}))))d\mu\\
&=\int_{\Omega}\Psi(\frac{f}{N_{\Phi}(f)})d\mu.
\end{align*}
Now for the case that $\mu(\Omega)=\infty$ and $\Psi\preceq \Phi$ globally, easily we get that
$$\int_{\Omega}\Psi(\frac{f}{N_{\Phi}(f)})d\mu\leq\int_{\Omega}\Phi(\frac{f}{N_{\Phi}(f)})d\mu\leq1.$$
And for the case that $\mu(\Omega)<\infty$ and $\Psi\preceq \Phi$, there exist $c>0$ and  $T>0$ such that for all $t\geq T$ we have $\Psi(t)\leq \Phi(t)$. Let $$E=\{w\in \Omega: \frac{f(w)}{N_{\Phi}(f)})\geq T\},\ \ \ \ \ \ \  N=\Psi(T)\mu(\Omega),$$ then we have
\begin{align*}
\int_{\Omega}\Psi(\frac{f}{N_{\Phi}(f)})d\mu&\leq \Psi(T)\mu(\Omega)+ \int_{\Omega\setminus E}\Phi(\frac{f}{N_{\Phi}(f)})d\mu\\
&\leq N+\int_{\Omega}\Phi(\frac{f}{N_{\Phi}(f)})d\mu\\
&\leq N+1.
\end{align*}
 So $N_{\Psi}(R_u(f))\leq CMN_{\Phi}(f)$, $N_{\Psi}(R_u(f))\leq CM(N+1)N_{\Phi}(f)$, respectively for infinite and finite cases. Thus $R_u$ is bounded in both cases and $\|R_u\|\leq C\|\Psi^{*^{-1}}(E(\Psi^*(u)))\|_{\infty}$, $\|R_u\|\leq C(N+1)\|\Psi^{*^{-1}}(E(\Psi^*(u)))\|_{\infty}$, respectively for infinite and finite cases.\\

\vspace*{0.3cm} {\bf Theorem 2.3.} Let WCT operator $R_{u}:\mathcal{D}\subseteq L^{\Phi}(\Sigma)\rightarrow L^{\Psi}(\Sigma)$ be well defined. Then the followings hold.\\

(a) Let $\Phi^*\circ\Psi^{*{-1}}\preceq \Theta$ globally for some Young's function $\Theta$. If $\Phi^*\in\bigtriangleup'$(globally), $\Theta\in\bigtriangledown'$(globally) and

(i) $E(\Phi^*(\bar{u}))=0$ on B,

(ii) $\sup_{n\in \mathbb{N}}\frac{E(\Phi^*(\bar{u}))(A_n)\mu(A_n)}{\Phi^*(\Psi^{*^{-1}}(\mu(A_n)))}<\infty$,

then $R_u$ is bounded. In another case, if $\Phi^*, \Psi^{*^{-1}}\in\bigtriangleup'$(globally),
 $$\sup_{n\in \mathbb{N}}E(\Phi^*(\bar{u}))(A_n)\mu(A_n)\Phi^*(\Psi^{*^{-1}}(\frac{1}{\mu(A_n)}))<\infty$$
 and (i) holds, then $R_u$ is bounded.\\

(b) If $\Phi^*\circ\Psi^{*^{-1}}$ is a Young's function, $\Phi^*\in \bigtriangledown'$ globally and $R_u$ is bounded from $L^{\Phi}(\Sigma)$ into $L^{\Psi}(\Sigma)$, then

(i) $E(\Phi^*(\bar{u}))=0$ on B,

(ii) $\sup_{n\in \mathbb{N}}E(\Phi^*(\bar{u}))(A_n)\mu(A_n)\Phi^*(\Psi^{*^{-1}}(\frac{1}{\mu(A_n)}))<\infty$.\\

{\bf Proof.} (a) If we prove that the operator $M_{\bar{u}}:L^{\Psi^*}(\mathcal{A})\rightarrow L^{\Phi^*}(\Sigma)$ is bounded, then we conclude that $R_u=(M_{\bar{u}})^*$ from $L^{\Phi}(\Sigma)$ into $L^{\Psi}(\Sigma)$ is bounded. So we prove the operator $M_{\bar{u}}$ is bounded under given conditions. Since $\Phi^*\in \bigtriangleup'$ and $\Theta\in\bigtriangledown'$ globally, then there exist $b, b'>0$ such that the following computations holds.
Put $M=\sup_{n\in \mathbb{N}}\frac{E(\Phi^*(\bar{u}))(A_n)\mu(A_n)}{\Phi^*(\Psi^{*^{-1}}(\mu(A_n)))}$. Therefore, for every $f\in L^{\Psi^*}(\mathcal{A})$ we get that
\begin{align*}
\int_{\Omega}\Phi^*(\frac{\bar{u}f}{N_{\Psi^*}(f)})d\mu&\leq b\int_{\Omega}E(\Phi^*(\bar{u}))\Phi^*(\frac{f}{N_{\Psi^*}(f)})d\mu\\
&=b\sum^{\infty}_{n=1}E(\Phi^*(\bar{u}))(A_n)\Phi^*(\frac{f}{N_{\Psi^*}(f)})(A_n)\mu(A_n)\\
&\leq b\sum^{\infty}_{n=1}E(\Phi^*(\bar{u}))(A_n)\Theta\Psi^*(\frac{f}{N_{\Psi^*}(f)})(A_n)\mu(A_n)\\
&\leq Mb\sum^{\infty}_{n=1}\Phi^*\Psi^{*^{-1}}(\mu(A_n))\Theta\Psi^*(\frac{f}{N_{\Psi^*}(f)})(A_n)\\
&\leq Mb\sum^{\infty}_{n=1}\Theta(\frac{\mu(A_n)}{b'})\Theta\Psi^*(\frac{f}{N_{\Psi^*}(f)})(A_n)\\
&\leq Mb\Theta(\frac{1}{b'})\int_{\Omega}\Psi^*(\frac{f}{N_{\Psi^*}(f)})d\mu\\
&\leq Mb \Theta(\frac{1}{b'}).
\end{align*}
This implies that $\int_{\Omega}\Phi^*(\frac{\bar{u}f}{N_{\Psi^*}(f)(Mb \Theta(\frac{1}{b'})+1)})d\mu\leq1$ and so $\| M_{\bar{u}}\|\leq (Mb \Theta(\frac{1}{b'})+1)$. Therefore the WCT operator $R_{u}$ is bounded. For the other case also by the same way we get that $R_u$ is bounded.\\

(b) Suppose that $R_u=M_{\bar{u}}^*$ from $L^{\Phi}(\Sigma)$ into $L^{\Psi}(\Sigma)$ is bounded, then the multiplication operator  $M_{\bar{u}}$ from $L^{\Psi^*}(\mathcal{A})$ into $L^{\Phi^*}(\Sigma)$ is  bounded. First, we show that $E(\Phi^*(\bar{u}))=0 $ on $B$.  Suppose on the contrary. Thus we can find some $\delta>0$ such that $\mu(\{w\in B:E(\Phi^*(\bar{u}))(w)>\delta\})>0$. Take $F=\{w\in
 B:E(\Phi^*(\bar{u}))(w)>\delta\}$. Since  $F\subseteq B$ is a $\mathcal{A}$-measurable set and $\mathcal{A}$ is $\sigma$-finite, then for each $n\in \mathcal{N}$, there exists $F_{n}\subseteq F$  with $F_n\in\mathcal{A}$ such that $\mu(F_{n})=\frac{\mu(F)}{2^n}$. Define $f_{n}=\frac{\chi_{F_{n}}}{N_{\Psi^*}(\chi_{F_{n}})}$. It is clear that $f_{n}\in L^{\Psi^*}(\mathcal{A})$ and $N_{\Psi^*}(f_n)=1$. Hence for each $k>0$ we have
\begin{align*}
\int_{\Omega}\Phi^*(\frac{k\bar{u}\chi_{F_{n}}}{N_{\Psi^*}(\chi_{F_{n}})})d\mu&=\int_{\Omega}\Phi^*(k\bar{u}\chi_{F_{n}}\Psi^{^*{-1}}(\frac{1}{\mu(F_{n})}))d\mu\\
&=\int_{\Omega}\Phi^*(k\bar{u}\Psi^{8^{-1}}(\frac{1}{\mu(F_{n})}))\chi_{F_{n}}d\mu\\
&\geq\int_{F_{n}}\Phi^*(\bar{u})\Phi^*(\frac{k\Psi^{*^{-1}}(\frac{1}{\mu(F_{n})})}{b})d\mu\\
&\geq\left(\int_{F_{n}}E(\Phi^*(\bar{u}))d\mu\right)\Phi^*(\frac{k}{b^2})\Phi^*(\Psi^{*^{-1}}(\frac{1}{\mu(F_{n})}))\\
&\geq \delta \Phi^*(\frac{k}{b^2})\Phi^*(\Psi^{*^{-1}}(\frac{1}{\mu(F_{n})}))\mu(F_n)\\
&=\delta \Phi^*(\frac{k}{b^2})\frac{\Phi^*(\Psi^{*^{-1}}(\frac{1}{\mu(F_{n})}))}{\frac{1}{\mu(F_{n})}}.
\end{align*}
Since $\Phi^*\circ\Psi^{*^{-1}}$ is a Young's function, then $\frac{\Phi^*(\Psi^{*^{-1}}(\frac{1}{\mu(F_{n})}))}{\frac{1}{\mu(F_{n})}}\rightarrow \infty$ when $n\rightarrow\infty$. Therefore $\int_{\Omega}\Phi^*(\frac{k\bar{u}\chi_{F_{n}}}{N_{\Psi^*}(\chi_{F_{n}})})d\mu\rightarrow \infty$ when $n\rightarrow\infty$ for each $k>0$ and so $N_{\Phi^*}(\bar{u}f_n)\rightarrow\infty$ when $n\rightarrow\infty$. This is a contradiction. It remains to prove (ii). Let $f_n=\frac{\chi_{A_{n}}}{N_{\Psi^*}(\chi_{A_{n}})}$, then we have
\begin{align*}
1\geq&\int_{\Omega}\Phi^*(\frac{\bar{u}f_n}{N_{\Phi^*}(\bar{u}f_n)})d\mu\\
&\geq\int_{\Omega}\Phi^*(\frac{\bar{u}\chi_{A_n}\Psi^{*^{-1}}(\frac{1}{\mu(A_n)})}{\|M_{\bar{u}}\|})d\mu\\
&\geq\int_{A_n}\Phi^*(\bar{u})\Phi^*(\frac{\Psi^{*^{-1}}(\frac{1}{\mu(A_n)})}{b\|M_{\bar{u}}\|})d\mu\\
&=\int_{A_n}E(\Phi^*(\bar{u}))\Phi^*(\frac{\Psi^{*^{-1}}(\frac{1}{\mu(A_n)})}{b\|M_{\bar{u}}\|})d\mu\\
&=E(\Phi^*(\bar{u}))(A_n)\Phi^*(\Psi^{*^{-1}}(\frac{1}{\mu(A_n)})\mu(A_n)\Phi^*(\frac{1}{b^2\|M_{\bar{u}}\|}).
\end{align*}
Hence we get that
\begin{align*}
\sup_{n}E(\Phi^*(\bar{u}))(A_n)\Phi^*(\Psi^{*^{-1}}(\frac{1}{\mu(A_n)})\mu(A_n)\leq\frac{1}{\Phi^*(\frac{1}{b^2\|M_{\bar{u}}\|})}<\infty.
\end{align*}
This completes the proof.\\

In the next proposition we give another necessary condition for boundedness of $R_u$. I think it's better that others.\\

{\bf Proposition 2.4.} Let WCT operator $R_{u}:\mathcal{D}\subseteq L^{\Phi}(\Sigma)\rightarrow L^{\Psi}(\Sigma)$ be well defined. And let $\Phi,\Psi\in\bigtriangleup'$ and $\Psi\circ\Phi^{-1}$ be a Young's function. The WCT operator $R_{u}$ into $L^{\Psi}(\Sigma)$ is bounded, if the following conditions hold;

(i) $E(\Phi^*(\bar{u}))=0$ on B,\\
(ii) $\sup_{n\in \mathbb{N}}\frac{\Phi\circ\Phi^{*^{-1}}(E(\Phi^*(u))(A_n))\mu(A_n)}{\Psi\circ\Phi^{-1}(\frac{1}{\mu(A_n)})}<\infty$.\\

In another case, if $\Psi\circ\Phi^{-1}$ isn't a Young's function, but $\Psi\circ\Phi^{-1}\preceq \Theta$ for some Young's function $\Theta$. Then the operator WCT operator $R_{u}$ into $L^{\Psi}(\Sigma)$ is bounded, if the conditions (i) and (ii) hold.\\
{\bf Proof.} Let $f\in L^{\phi}(\Omega)$ such that $N_{\Phi}(f)\leq1$. Then we have
\begin{align*}
\int_{\Omega}\Psi(\frac{E(uf)}{N_{\Phi}(f)})d\mu&\leq\int_{\Omega}\Psi\left(\Phi^{-1}(E(\Phi(\frac{f}{N_{\Phi}(f)})))\Phi^{*^{-1}}(E(\Phi^*(u))))\right)d\mu\\
&=\sum^{\infty}_{n=1}\Psi\left(\Phi^{-1}(E(\Phi(\frac{f}{N_{\Phi(f)}})))(A_n)\right)\Psi\left(\Phi^{*^{-1}}(E(\Phi^*(u)))(A_n))\right)\mu(A_n)\\
&\leq M\sum^{\infty}_{n=1}\Psi\circ\Phi^{-1}\left((E(\Phi(\frac{f}{N_{\Phi}(f)})))(A_n)\mu(A_n)\right)\\
&\leq M\sum^{\infty}_{n=1}\Psi\circ\Phi^{-1}\left((E(\Phi(\frac{f}{N_{\Phi}(f)})))(A_n)\mu(A_n)\right)\\
&\leq M\Psi\circ\Phi^{-1}\left(\int_{\Omega}E(\Phi(\frac{f}{N_{\Phi}(f)}))d\mu\right)\\
&=M\Psi\circ\Phi^{-1}\left(\int_{\Omega}\Phi(\frac{f}{N_{\Phi}(f)})d\mu\right)\\
&\leq M\Psi\circ\Phi^{-1}(1).
\end{align*}
Hence $$\int_{\Omega}\Psi(\frac{E(uf)}{(M\Psi\circ\Phi^{-1}(1)+1)N_{\Phi}(f)})d\mu\leq1.$$
Consequently we get that $$N_{\Psi}(E(uf))\leq (M\Psi\circ\Phi^{-1}(1)+1)N_{\Phi}(f).$$
 Thus the operator $R_u$ is bounded. If $\Psi\circ\Phi^{-1}\preceq \Theta$ for some Young's function $\Theta$, by the same method we get that $$N_{\Psi}(E(uf))\leq (M\Theta(1)+1)N_{\Phi}(f).$$ So this also states that the operator $R_u$ is bounded.\\


 {\bf Remark 2.5.} If $(\Omega, \Sigma, \mu)$ is a non-atomic measure space, then under assumptions of Theorem 2.3, there is not any non-zero bounded operator of the form $R_u$ from $L^{\Phi}(\Sigma)$ into $L^{\Psi}(\Sigma)$.\\

Put $\Phi(x)=\frac{x^p}{p}$ for $x\geq0$, where $1<p<\infty$. It is clear that $\Phi$ is a Young's function and $\Phi^*(x)=\frac{x^{p'}}{p'}$, where $1<p'<\infty$ and $\frac{1}{p}+\frac{1}{p'}=1$. If $\Psi(x)=\frac{x^q}{q}$ for $x\geq0$, where $1<p<q<\infty$. Then $\Phi^*\circ\Psi^{*^{-1}}(x)=\frac{q'^{\frac{1}{q'}}}{p'}x^{\frac{p'}{q'}}$. Since $\frac{p'}{q'}>1$, then $\Phi^*\circ\Psi^{*^{-1}}$ is a Young's function. These observations and Theorems 2.3, 2.4 give us the next Remark.\\

{\bf Remark 2.6.} Let $R_{u}:\mathcal{D}\subseteq L^{p}(\Sigma)\rightarrow L^{q}(\Sigma)$ be well defined. Then the operator $R_{u}$ from $L^{p}(\Sigma)$ into $L^{q}(\Sigma)$, where $1<p<q<\infty$, is bounded if and only if the followings hold:

(i) $E(|u|^{p'})=0$ on $B$.

(ii) $\sup_{n\geq0}\frac{E(|u|^{p'})(A_n)}{\mu(A_n)^{\frac{p'}{q'}-1}}<\infty$.\\
In addition, we get that there is not any non-zero bounded operator of the form $R_u$ from $L^p$ into $L^q$, ($1<p<q<\infty$) when the underlying measure space is non-atomic.

Specially, if $\mathcal{A}=\Sigma$, then $E=I$ and so the multiplication operator $M_u$ from $L^{p}(\Sigma)$ into $L^{q}(\Sigma)$ is bounded if and only if

(i) $u=0$ on $B$.

(ii) $\sup_{n\geq0}\frac{u(A_n)}{\mu(A_n)^{\frac{1}{q'}-\frac{1}{p'}}}<\infty$.\\

Here we recall a fundamental Lemma, which is as an easy exercise.

{\bf Lemma 2.7.} Let $\Phi_i$, $i=1,2,3$, be Young's functions for which
$$\Phi_3(xy)\leq \Phi_1(x)+\Phi_2(y), \ \ \ \ \ \ x\geq0, y\geq0.$$
If $f_i\in L^{\Phi_i}(\Sigma)$, $i=1,2$, where $(\Omega,\Sigma,\mu)$ is any measure space, then
$$N_{\Phi_3}(f_1f_2)\leq2N_{\Phi_1}(f_1)N_{\Phi_2}(f_2).$$

{\bf Theorem 2.8.} Let $\Phi$ and $\Psi$ be Young's functions an $R_{u}:\mathcal{D}\subseteq L^{\Phi}(\Sigma)\rightarrow L^{\Psi}(\Sigma)$ be well defined. Then the followings hold:\\

 (i) Suppose that there exists a Young's function $\Theta$ such that
 $$\Psi(xy)\leq \Phi(x)+\Theta(y), \ \ \ \ \ \ x\geq0, y\geq0$$
 and $(E,\Phi)$ satisfies GCH-inequality. In this case if $\Phi^{*^{-1}}(E(\Phi^*(u)))\in L^{\Theta}(\mathcal{A})$, then the WCT operator $R_u$ from $L^{\Phi}(\Sigma)$ into $L^{\Psi}(\Sigma)$ is bounded.\\

 (ii) Let $\Theta=\Psi^*\circ\Phi^{*^{-1}}$ be a Young's function, $\Theta\in \bigtriangleup_2$ and $\Phi^* \in \bigtriangleup_2$. In this case if WCT operator $R_u$ is bounded from $L^{\Phi}(\Sigma)$ into $L^{\Psi}(\Sigma)$, then $E(\Phi^*(\bar{u}))\in L^{\Theta^*}(\mathcal{A})$. Consequently $\Phi^{*^{-1}}(E(\Phi^*(\bar{u})))\in L^{\Theta^*\circ\Phi^*}(\mathcal{A})$.\\

 {\bf Proof.} (i) Let $f\in L^{\Phi}(\Sigma)$ such that $N_{\Phi}(f)\leq1$. This means that
 $$\int_{\Omega}\Phi(\Phi^{-1}(E(\Phi(f))))d\mu=\int_{\Omega}\Phi(f)d\mu\leq1,$$
 hence $N_{\Phi}(\Phi^{-1}(E(\Phi(f))))\leq1$. By using GCH-inequality we have
 \begin{align*}
 N_{\Psi}(E(uf))&\leq N_{\Phi}(\Phi^{-1}(E(\Phi(f))))N_{\Theta}(\Phi^{*^{-1}}(E(\Phi^*(u))))\\
 &\leq N_{\Theta}(\Phi^{*^{-1}}(E(\Phi^*(u)))).
 \end{align*}
Thus for all $f\in L^{\Phi}(\Sigma)$ we have $$N_{\Psi}(E(uf))\leq N_{\Phi}(f)N_{\Theta}(\Phi{*^{-1}}E(\Phi^*(u)).$$ And so the operator $R_u$ is bounded.\\

(ii) Suppose that $R_u$ is bounded. So the adjoint operator $M_{\bar{u}}=(R_u)^{\ast}:L^{\Psi^*}(\mathcal{A})\rightarrow L^{\Phi^*}(\Sigma)$ is also bounded. For $f\in L^{\Theta}(\mathcal{A})$ we have $\Phi^{*^{-1}}(f)\in L^{\Psi^*}(\mathcal{A})$. Consequently we get that
\begin{align*}
\int_{\Omega}E(\Phi^*(\bar{u}))fd\mu&=\int_{\Omega}\Phi^*(\bar{u})\Phi^*(\Phi^{*^{-1}}(f))d\mu\\
&\leq b\int_{\Omega}\Phi^*(\bar{u}\Phi^{*^{-1}}(f))d\mu\\
&=b\int_{\Omega}\Phi^*(M_{\bar{u}}(\Phi^{*^{-1}}(f)))d\mu<\infty.
\end{align*}
Therefore $\int_{\Omega}E(\Phi^*(\bar{u}))fd\mu<\infty$ for all $f\in L^{\Theta}(\mathcal{A})$. This implies that $E(\Phi^*(\bar{u}))\in L^{\Theta^*}(\mathcal{A})$. This completes the proof.\\

\vspace*{0.3cm} {\bf Remark 2.9.} Let $R_{u}:\mathcal{D}\subseteq L^{p}(\Sigma)\rightarrow L^{q}(\Sigma)$ be well defined. Then the operator $R_{u}$ from $L^{p}(\Sigma)$ into $L^{q}(\Sigma)$, where $1<q<p<\infty$, is bounded if and only if $(E(|u|^{p'}))^{\frac{1}{p'}}\in L^{r}(\mathcal{A})$, where $r=\frac{pq}{p-q}$.

Specially, if $\mathcal{A}=\Sigma$, then $E=I$ and so the multiplication operator $M_u$ from $L^{p}(\Sigma)$ into $L^{q}(\Sigma)$ is bounded if and only if $u\in L^{r}(\Sigma)$.\\

\vspace*{0.3cm} {\bf Example 2.10.} Let $\Omega=[-1,1]$, $d\mu=\frac{1}{2}dw$ and
$\mathcal{A}=\langle\{(-a,a):0\leq a\leq1\}\rangle$
($\sigma$-algebra generated by symmetric intervals). Then
 $$E^{\mathcal{A}}(f)(w)=\frac{f(w)+f(-w)}{2}, \ \ w\in \Omega,$$
 where $E^{\mathcal{A}}(f)$ is defined. Thus $E^{\mathcal{A}}(|f|)\geq\frac{|f|}{2}$. Hence $|f|\leq2E(|f|)$.
 Let $\Phi(w)=e^{w^p}-w^p-1$ and $\Psi(w)=\frac{w^p}{p}$ be Young's functions, where $p>1$. For each $f\in L^{\Phi}(\Omega, \Sigma, \mu)$ we have  $\Phi(|f|)\leq2E(\Phi(|f|))$.
 This implies that
 $$E(|fg|)\leq 4 \Phi^{-1}(E(\Phi(|f|)))\Psi^{-1}(E(\Psi(|g|))).$$
 If $u$ is a non-zero continuous function on $\Omega$,
 then for Young's function $\Theta(w)=(1+w^p)log(1+w^p)-w^p$ we have
 $$\Psi(xy)\leq \Phi(x)+\Theta(y), \ \ \ \ \ \ -1\leq x,y\leq1.$$
So by Theorem 2.8 the WCT operator $R_u$ is bounded from $L^{\Phi}$ into $L^{\Psi}$. But it is not bounded from $L^{\Psi}$ into $L^{\Phi}$, because of Theorem 2.3.\\
\vspace*{0.3cm} {\bf Example 2.11.}
Let $\Omega=[0,1]$, $\Sigma$ be the $\sigma$-algebra of
Lebesgue measurable subset of $\Omega$ and let $\mu$ be the Lebesgue
measure on $\Omega$. Fix $n\in\{2, 3, 4 . . .\}$ and let $s:[0,
1]\rightarrow[0, 1]$ be defined by $s(w)=w+\frac{1}{n}$(mod 1).
Let $\mathcal{B}=\{E\in \Sigma:s^{-1}(E)=E\}$. In this case
$$E^{\mathcal{B}}(f)(w)=\sum^{n-1}_{j=0}f(s^{j}(w)),$$ where $s^j$
denotes the jth iteration of $s$. The functions $f$ in the range
of $E^{\mathcal{B}}$ are those for which the n graphs of f
restricted to the intervals $[\frac{j-1}{n},\frac{j}{n}]$, $1\leq
j\leq n$, are all congruent. If $\Phi(w)=e^{w^4}-1$ and $\Psi(w)=\frac{w^2}{log(e+w)}$, then for $\Theta(w)=\Phi^{*}(w^2)$ we have

$$\Psi(xy)\leq \Phi(x)+\Theta(y), \ \ \ \ \ \ 0\leq x,y\leq1.$$
Then by Theorem 2.8 $R_u$ is a bounded operator from $L^{\Phi}$ into $L^{\Psi}$ for every non-zero continuous function $u$. But it is not bounded from $L^{\Psi}$ into $L^{\Phi}$, because of Theorem 2.3.\\


\section{ \sc\bf Essential norm }
Let $\mathfrak{B}$ be a Banach space and $\mathcal{K}$ be the set
of all compact operators on $\mathfrak{B}$. For $T\in
L(\mathfrak{B})$, the Banach algebra of all bounded linear
operators on $\mathfrak{B}$ into itself, the essential norm of $T$
means the distance from $T$ to $\mathcal{K}$ in the operator norm,
namely $\|T\|_e =\inf\{\|T - S\| : S \in\mathcal{K}\}$. Clearly,
$T$ is compact if and only if $\|T\|_e= 0$. Let $X$ and $Y$ be
reflexive Banach spaces and $T\in L(X,Y)$. It is easy to see that
$\|T\|_{e}=\|T^*\|_{e}$. In this section we assume that $a_j=\mu(A_j)$, where
$A_j$'s are $\mathcal{A}$-atoms.\\

 In the sequel we present an upper bound for
essential norm of $EM_u$ on Orlicz space $L^{\Phi}(\Sigma)$. For this we first recall some results of \cite{e3} for compactness of WCT operator $R_u$.\\

{\bf Theorem 3.1.}\cite{e3} Let  WCT operator $R_{u}$ be bounded on
$L^{\Phi}(\Sigma)$, then the followings hold.\\

(a) If $R_u$ is compact, then
 $$N_{\varepsilon}(E(u))=\{w\in \Omega:E(u)(w)\geq\varepsilon\}$$
 consists of finitely many $\mathcal{A}-$atoms, for all  $\varepsilon>0$.\\

(b) If $R_u$ is compact and $\Phi\in \bigtriangleup'$(globally), then $N_{\varepsilon}(\Phi^{*^{-1}}(E(\Phi^*(u))))$

consists of finitely many $\mathcal{A}-$atoms, for all $\varepsilon>0$, where

$$N_{\varepsilon}(\Phi^{*^{-1}}(E(\Phi^*(u))))=\{w\in \Omega:\Phi^{*^{-1}}(E(\Phi^*(u)))(w)\geq\varepsilon\}.$$\\

(c) If $(E, \Phi)$ satisfies the GCH-inequality and $N_{\varepsilon}(\Phi^{*^{-1}}(E(\Phi^*(u))))$
consists of finitely many $\mathcal{A}-$atoms, for all
$\varepsilon>0$, then $T$ is compact.\\

{\bf Corollary 3.2.}\cite{e3} Under assumptions of theorem 3.1 we have the followings:\\
(a) If  $(E, \Phi)$ satisfies the
GCH-inequality and $\Phi\in
\bigtriangleup'$(globally), then $T$ is compact if and only if
 $N_{\varepsilon}(\Phi^{*^{-1}}(E(\Phi^*(u))))$ consists of finitely many $\mathcal{A}-$atoms, for all
 $\varepsilon>0$.\\

(b) If $\Phi^*\prec x$(globally) and $(E, \Phi)$ satisfies the
GCH-inequality, then then $T$
is compact if and only if
 $N_{\varepsilon}(\Phi^{*^{-1}}(E(\Phi^*(u))))$ consists of finitely many $\mathcal{A}-$atoms, for all
 $\varepsilon>0$.\\

(c)If $(\Omega, \mathcal{A}, \mu)$ is non-atomic measure
space,$(E, \Phi)$ satisfies the GCH-inequality and $\Phi\in \bigtriangleup'$(globally).
Then $R_{u}$ is a compact operator on $L^{\Phi}(\Sigma)$ if and
only if $R_u=0$.\\

\vspace*{0.3cm} {\bf Theorem 3.3.} Let $R_{u}: L^{\Phi}(\Sigma)\rightarrow L^{\Phi}(\Sigma)$ is bounded. Then\\

(a) If $(E, \Phi)$ satisfies in GCH-inequality and  $\beta_2=\inf\{\varepsilon>0:N_{\varepsilon}$
consists of finitely many $\mathcal{A}$-atoms$\}$, where
$N_{\varepsilon}=N_{\varepsilon}(\Phi^{*^{-1}}(E(\Phi^*(u))))$. Then

$$\|R_{u}\|_{e}\leq C\beta_2,$$
where $C$ comes from GCH-inequality.\\

(b) If $a_n\rightarrow 0$ or $\{a_n\}_{n\in \mathbb{N}}$
has no convergent subsequence. Let $\beta_1=\inf\{\varepsilon>0:N_{\varepsilon}$
consists of finitely many $\mathcal{A}$-atoms$\}$, where
$N_{\varepsilon}=N_{\varepsilon}(E(u))$. Then

$$\|R_{u}\|_{e}\geq \beta_1.$$

{\bf Proof} (a) Let $\varepsilon>0$. Then $N_{\varepsilon+\beta_2}$
consist of finitely many $\mathcal{A}$-atoms. Put
$u_{\varepsilon+\beta_2}=u\chi_{N_{\varepsilon+\beta_2}}$ and
$R_{u_{\varepsilon+\beta_2}}$. So $R_{u_{\varepsilon+\beta_2}}$ is
finite rank and so compact. And for every $f\in L^{\Phi}(\Sigma)$ we have\\
\begin{align*}
\int_{\Omega}\Phi(\frac{R_u(f)-R_{\varepsilon+\beta_2}(f)}{C(\varepsilon+\beta_2) N_{\Phi}(f)})d\mu&=\int_{\Omega}\Phi(\frac{E(uf)\chi_{\Omega \setminus N_{\varepsilon+\beta}}}{C(\varepsilon+\beta_2) N_{\Phi}(f)})d\mu\\
&\leq\int_{\Omega\setminus N_{\varepsilon+\beta_2}}\Phi(\frac{C \Phi^{-1}(E(\Phi(|\frac{f}{N_{\Phi}(f)})|)))\Phi^{*^{-1}}(E(\Phi^*(|u|)))}{C(\varepsilon+\beta_2)})d\mu\\
&\leq\int_{\Omega\setminus N_{\varepsilon+\beta_2}}\Phi(\Phi^{-1}(E(\Phi(|\frac{f}{N_{\Phi}(f)})|)))d\mu\leq\int_{\Omega}E(\Phi(\frac{f}{N_{\Phi}(f)}))d\mu\\
&=\int_{\Omega}\Phi(\frac{f}{N_{\Phi}(f)})d\mu\leq1.
\end{align*}
This implies that

$$\|R_u\|_{e}\leq \|R_u-R_{u_{\varepsilon+\beta_2}}\|\leq C(\beta_2+\varepsilon).$$

This mean's that $\|R_{u}\|_{e}\leq C\beta_2$.\\

(b)  Let $0<\varepsilon<\beta_1$. Then by definition,
$N_{\beta_1-\varepsilon}=N_{\beta_1-\varepsilon}(E(u))$ contains infinitely many atoms or a
non- atomic subset of positive measure. If
$N_{\beta_1-\varepsilon}$ consists a non- atomic subset, then
we can find a sequence $\{B_{n}\}_{n\in \mathbb{N}}$ such that
$\mu(B_{n})<\infty$ and $\mu(B_{n})\rightarrow 0$. Put
$f_{n}=\frac{\chi_{B_{n}}}{N_{\Phi}(\chi_{B_{n}})}$, then for
every $A\in \Sigma$ with $0<\mu(A)<\infty$ we have

$$\int_{\Omega}f_{n}\chi_{A}d\mu=\mu(A\cap B_{n})\Phi^{-1}(\frac{1}{\mu(B_{n})})\leq \frac{\Phi^{-1}(\frac{1}{\mu(B_{n})})}{\frac{1}{\mu(B_{n})}}
\rightarrow0.$$ when $n\rightarrow \infty$. Also, if
$N_{\beta_1-\varepsilon}$ consists infinitely many atoms
$\{A'_{n}\}_{n\in \mathbb{N}}$. We set
$f_{n}=\frac{\chi_{A'_{n}}}{N_{\Phi}(\chi_{A'_{n}})}$. Then
for every $A\in \Sigma$ with $0<\mu(A)<\infty$ we have
$$\int_{\Omega}f_{n}\chi_{A}d\mu=\mu(A\cap A'_{n})\Phi^{-1}(\frac{1}{\mu(A'_{n})}).$$
If $\{\mu(A_{n})\}_{n\in \mathbb{N}}$ has no convergent
subsequence, then there exists $n_{0}$ such that for $n>n_{0}$,
$\mu(A\cap A_{n}')=0$ and if $\mu(A_{n})\rightarrow 0$ then
$\mu(A_{n}')\rightarrow 0$. Thus
$\int_{\Omega}f_{n}\chi_{A}d\mu=\mu(A\cap
A'_{n})\Phi^{-1}(\frac{1}{\mu(A'_{n})})\rightarrow0$ in both
cases. These imply that $f_{n}\rightarrow0$ weakly. So
$$\int_{\Omega}\Phi(\frac{(\beta_1-\varepsilon)f_{n}}{N_{\Phi}(R_u(f_{n}))})d\mu
\leq\int_{\Omega}\Phi(\frac{E(u)f_{n}}{N_{\Phi}(R_u(f_{n}))})d\mu=\int_{\Omega}\Phi(\frac{R_u(f_{n})}{N_{\Phi}(R_u(f_{n}))})d\mu.$$

Thus $N_{\Phi}(R_u(f_{n}))\geq \beta_2-\varepsilon$.

Also, there exists compact operator $T\in L(L^{\Phi}(\Sigma))$ such that
$\|R_{u}\|_{e}\geq\|T-R_{u}\|-\varepsilon$. Hence
$N_{\Phi}(Tf_{n})\rightarrow o$ and so there exists $N>0$ such
that for each $n>N$, $N_{\Phi}(Tf_{n})\leq \varepsilon$. So

$$\|R_{u}\|_{e}\geq\|R_{u}-T\|-\varepsilon\geq|N_{\Phi}(R_u(f_{n}))-N_{\Phi}(Tf_{n})|\geq\beta_1-\varepsilon-\varepsilon,$$
thus we conclude that $\|R_{u}\|_{e}\geq\beta_1$.\\

\vspace*{0.3cm} {\bf Corollary 3.4.} Let
$u:\Omega\rightarrow\mathbb{C}$ be $\Sigma-$measurable and Let
$M_{u}: L^{\Phi}(\Sigma)\rightarrow
L^{\Phi}(\Sigma)$. If
$\beta=\beta_1=\beta_2=\inf\{\varepsilon>0:N_{\varepsilon}$ consists of finitely
many atoms$\}$. Then\\

(a) $\|M_{u}\|_{e}\leq \beta$.\\

(b) Let $\Phi\in \bigtriangleup_{2}$ and
$a_n\rightarrow 0$ or $\{a_n\}_{n\in \mathbb{N}}$
has no convergent subsequence. Then $\beta=\|M_{u}\|_{e}$.\\


\begin{thebibliography}{99}




\bibitem{besh} C. Bennett and R. Sharpley, Interpolation of
operators, Academic Press, INC, 1988.



\bibitem{e3} Y. Estaremi, Multiplication conditional expectation type operators on Orlicz spaces, J. Math. Anal. Appl. {\bf 414} (2014) 88–98.



%
%
%
%
%



\bibitem{rao}
M. M. Rao, Conditional measure and applications, Marcel Dekker,
New York, 1993.

\bibitem{raor} M.M. Rao, Z.D. Ren, Theory of Orlicz spaces, Marcel Dekker,
New York, 1991.

%
%
%
%
%
%
%

\bibitem{taky}
H. Takagi and K. Yokouchi, Multiplication and composition operators
between two $L^p$-spaces, Contemporary Math. {\bf 232}(1999),
321-338.

\bibitem{z}
A. C. Zaanen, Integration, 2nd ed., North-Holland, Amsterdam,
1967.
%

\end{thebibliography}
\end{document}